\documentclass[onecolumn,amsmath,amssymb,preprint]{revtex4}
\usepackage{mathtext,bm,bbm,amsmath,amsfonts,amssymb,indentfirst,syntonly,graphicx}
\usepackage{graphicx}
\usepackage{epsfig}
\usepackage{graphicx}

\begin{document}

\title{Finslerian Ressiner-Nordstrom spacetime}

\author{Xin Li $^{1}$}
\email{lixin1981@cqu.edu.cn}
\affiliation{$^1$Department of Physics, Chongqing University, Chongqing 401331, China}

\begin{abstract}
We have obtained Finslerian Ressiner-Nordstrom solution where it is asymptotic to a Finsler spacetime with constant flag curvature while $r\rightarrow\infty$. The covariant derivative of modified Einstein tensor in Finslerian gravitational field equation for this solution is conserved. The symmetry of the special Finslerian Ressiner-Nordstrom spacetime, namely, Finsler spacetime with constant flag curvature, has been investigated. It admits four independent Killing vectors. The Finslerian Ressiner-Nordstrom solution differs from Ressiner-Nordstrom metric only in two dimensional subspace. And our solution requires that its two dimensional subspace have constant flag curvature. We have obtained eigenfunction of Finslerian Laplace operator of ``Finslerian sphere", namely, a special subspace with positive constant flag curvature. The eigenfunction is of the form $\bar{Y}_l^m=Y_l^m+\epsilon^2(C_{l+2}^m Y_{l+2}^m+C_{l-2}^m Y_{l-2}^m)$ in powers of Finslerian parameter $\epsilon$, where $C_{l+2}^m$ and $C_{l-2}^m$ are constant. However, the eigenvalue depends on both $l$ and $m$. The eigenvalues corresponded to $Y_1^0$ remain the same with Riemannian Laplace operator and the eigenvalues corresponded to $Y_1^{\pm1}$ are different. This fact just reflect the symmetry of ``Finslerian sphere", which admits a $z$-axis rotational symmetry and breaks other symmetry of Riemannian sphere. The eigenfunction of Finslerian Laplace operator implies that monopolar and dipolar terms of multipole expansion of gravitational potential are unchanged and other multipole terms are changed.
\end{abstract}
%\pacs{04.30.-w, 98.65.-r, 98.80.-k}
%\keywords{gravitational wave; standard siren; anisotropy of the universe}

\maketitle

\section{Introduction}
A black hole is a specific region of spacetime that has strong enough gravity such that even the light cannot escape from it. The black hole physics has been discussed intensively by physicists. The researches on black hole physics origins from the exact solution of Einstein's gravitational field equation. In four dimensional spacetime, Schwarzschild and Kerr solutions are exact solutions of Einstein vacuum field equation, which correspond to the spherical symmetry and axis symmetry respectively. And Reissner-Nordstrom spacetime is a solution that corresponds to the gravitational field generated by a charged and spherical symmetric gravitational source. General form of the solutions are called Kerr-Newmann spacetime, which corresponds to the gravitational field generated by a charged and axis symmetric gravitational source \cite{Kerr-Newmann}. Inspiring from Kerr-Newmann spacetime, the no hair theorem of the black hole is found by physicists \cite{Gravitation}. It states that the black hole only has three properties, namely, mass, electric charge and spin. The above solutions or spacetimes are asymptotically flat. Schwarzschild-de Sitter spacetime is a solution with cosmological horizons, and it is asymptotic to de Sitter spacetime \cite{Bousso}. Based on the Kerr-Newmann spactime, Hawking and Bekenstein et al. involved the concept of entropy for black hole \cite{Benkenstein} and constructed theory of black hole thermodynamics, then proposed four laws of black hole thermodynamics \cite{Bardeen}. Hawking radiation \cite{Hawking radiation}, as a important complement of black hole thermodynamics, exhibits quantum feature of black hole.

Three physical processes could lead to formations of black hole. One is the gravitational collapse of a heavy star\cite{Celotti}. The second is the gravitational collapse of a primordial overdensity in the early universe\cite{Carr1,Carr2,Carr3,Clesse,Dolgov}. The third is high-energy collisions \cite{Giddings}. The merge of black hole binary could generate gravitational wave, such phenomena has been observed by the Advanced LIGO detectors, i.e. GW150914, GW151226, GW170104 and GW170814 \cite{GW150914,GW151226,GW170104,GW170814}.

It is interesting and important to search more exact solutions of gravitational field equation in four dimensional spacetime. And the black holes that corresponded to these solutions are expected to be tested by the gravitational wave detectors in the near future. Finsler geometry \cite{Book by Bao} is a natural generation of Riemannian geometry. The basic feature of Finsler geometry is that its length element does not have quadratic restriction. Generally, Finslerian extension of a given Riemannian spacetime has lower symmetry than the Riemannian spacetime\cite{Wang,Finsler PF}. A typical example of Finsler spacetime, i.e., Randers spacetime \cite{Randers}, breaks rotational symmetry and induces parity violation. By this basic feature of Finsler geometry, Finsler geometry is used to describe violation of Lorentz invariance \cite{Girelli,Gibbons,Finsler LI,Kostelecky} and study the anisotropy of our universe \cite{Kouretsis,Finsler alphae}.

We have suggested an anisotropic inflation model in which the background space is taken to be a Randers space. This anisotropic inflation model could account for power asymmetry of the cosmic microwave background (CMB)\cite{Finsler inflation}. In Finsler geometry, there is no unique extension of Riemannian geometric objects, such as the connections and curvature\cite{Book by Bao}. This feature also involves in searching the gravitational field equation in Finsler spacetime \cite{Li Berwald,Miron,Rutz,Vacaru,Pfeifer Gravity}. We have proposed the gravitational field equation in Finsler spacetime and found a non-Riemannian exact solution \cite{Finsler BH}. This solution or metric is none other than the Schwarzschild metric except for the change from the Riemann sphere to ``Finslerian sphere". And interior solution for the Finslerian Schwarzschild metric also exists. We have proved that covariant derivative of the Finslerian gravitational field equation for the metric is conserved. It is interesting to search a general solution for our Finslerian gravitational field equation and test if black hole corresponded to the solution possess the three properties, namely, mass, electric charge and spin or not. The spherical harmonics are eigenfunction of Laplace operator for the Riemann sphere, which play an important role in modern physics. Since our Finslerian Schwarzschild solution admits a ``Finslerian sphere", it is worth investigating eigenfunction of Laplace operator for ``Finslerian sphere". It is expected that symmetry of ``Finslerian sphere" has direct influence on its eigenfunction.

This paper is organized as follows. In section II, we first give a brief introduction to Ressiner-Nordstrom metric. Then, we present Finslerian Ressiner-Nordstrom solution in Finsler spacetime. We discuss the symmetry of Finslerian Ressiner-Nordstrom metric at the last of section II. In sec III, we present Finslerian Laplace operator of ``Finslerian sphere", and give the eigenfunction and corresponded eigenvalue of Finslerian Laplace operator. Conclusions and remarks are given in section IV.

\section{Exact solution of gravitational field equation in Finsler spacetime}
\subsection{Brief introduction to Ressiner-Nordstrom metric}
In general relativity, the Ressiner-Nordstrom spacetime is given as follows
\begin{equation}
ds^2=-fdt^2+f^{-1}dr^2+r^2d\Omega^2_k,
\end{equation}
where $f=k-\frac{2GM}{r}-br^2+\frac{4\pi GQ^2}{r^2}$, and the metric $d\Omega^2_k$ denotes the 2 dimensional metric with constant sectional curvature $k$. Usually, after a reparameterization, we can set $k$ to be $1,0,-1$.
The Ressiner-Nordstrom is a solution of Einstein field equation, i.e.,
\begin{equation}
Ric_{\mu\nu}-g_{\mu\nu}S/2=8\pi G T_{\mu\nu},
\end{equation}
where $Ric_{\mu\nu}$ is Ricci tensor and $S=g^{\mu\nu}Ric_{\mu\nu}$ is the scalar curvature. The energy momentum tensor of the Ressiner-Nordstrom spacetime is given as
\begin{equation}
T_{\mu\nu}=T^{em}_{\mu\nu}+T^{c}_{\mu\nu}
\end{equation}
where $T^{em}_{\mu\nu}=\frac{Q^2}{2r^4}diag\{f,-f,r^2g^\omega_{ij}\}$ ($g^\omega_{ij}$ is the metric of $d\Omega^2_k$)denotes the energy momentum tensor of electromagnetic field and $T^{c}_{\mu\nu}=-3b g_{\mu\nu}/8\pi G$. The Ressiner-Nordstrom spacetime will reduce to four dimensional spacetime with constant curvature, namely, the de Sitter spacetime, if $Q=M=0$ and $k=1$.

\subsection{Finslerian Ressiner-Nordstrom solution}
Instead of defining an inner product structure over the tangent bundle in Riemann geometry, Finsler geometry is based on
the so called Finsler structure $F$ with the property
$F(x,\lambda y)=\lambda F(x,y)$ for all $\lambda>0$, where $x\in M$ represents position
and $y\equiv\frac{dx}{d\tau}$ represents velocity. The Finslerian metric is given as\cite{Book
by Bao}
\begin{equation}\label{Finsler metric}
g_{\mu\nu}\equiv\frac{\partial}{\partial
y^\mu}\frac{\partial}{\partial y^\nu}\left(\frac{1}{2}F^2\right).
\end{equation}

In our previous research \cite{Finsler BH}, we have proposed that the vacuum gravitational field equation is vanish of the Ricci scalar, and we have obtained a solution of Finslerian vacuum field equation. It is given as
\begin{equation}
F^2=-\left(1-\frac{2GM}{r}\right)y^ty^t+\left(1-\frac{2GM}{r}\right)^{-1}y^r y^r+r^2F^2_{FS},
\end{equation}
where $F_{FS}$ is a two dimensional Finsler spacetime with positive constant flag curvature. The Ricci scalar is given as
\begin{equation}\label{Ricci scalar}
Ric\equiv R^\mu_{~\mu}=\frac{1}{F^2}\left(2\frac{\partial G^\mu}{\partial x^\mu}-y^\lambda\frac{\partial^2 G^\mu}{\partial x^\lambda\partial y^\mu}+2G^\lambda\frac{\partial^2 G^\mu}{\partial y^\lambda\partial y^\mu}-\frac{\partial G^\mu}{\partial y^\lambda}\frac{\partial G^\lambda}{\partial y^\mu}\right),
\end{equation}
where
\begin{equation}
\label{geodesic spray}
G^\mu=\frac{1}{4}g^{\mu\nu}\left(\frac{\partial^2 F^2}{\partial x^\lambda \partial y^\nu}y^\lambda-\frac{\partial F^2}{\partial x^\nu}\right)
\end{equation} is called geodesic spray coefficients.

Now, we propose an ansatz that the Finsler structure is of the form
\begin{equation}\label{Schwarz like}
F^2=-f(r)y^ty^t+f(r)^{-1}y^r y^r+r^2\bar{F}^2(\theta,\varphi,y^\theta,y^\varphi).
\end{equation}
Throughout the paper, the index labeled by Greek alphabet denote the index of four dimensional spacetime $F$, and the index labeled by Latin alphabet denote the index of two dimensional subspace $\bar{F}$. Plugging the Finsler structure (\ref{Schwarz like}) into the formula (\ref{geodesic spray}), we obtain that
\begin{eqnarray}\label{geodesic spray t}
G^t&=&\frac{f'}{2f}y^t y^r,\\
\label{geodesic spray r}
G^r&=&-\frac{f'}{4f}y^r y^r+\frac{ff'}{4}y^t y^t-\frac{r}{2A}\bar{F}^2,\\
\label{geodesic spray theta}
G^\theta&=&\frac{1}{r}y^\theta y^r+\bar{G}^\theta,\\
\label{geodesic spray phi}
G^\varphi&=&\frac{1}{r}y^\varphi y^r+\bar{G}^\varphi,
\end{eqnarray}
where the prime denotes the derivative with respect to $r$, and the $\bar{G}$ is the geodesic spray coefficients derived by $\bar{F}$.
Plugging the geodesic coefficients (\ref{geodesic spray t},\ref{geodesic spray r},\ref{geodesic spray theta},\ref{geodesic spray phi}) into the formula of Ricci scalar (\ref{Ricci scalar}), we obtain that
\begin{eqnarray}\label{Ricci scalar1}
F^2Ric=\left[\frac{ff''}{2}+\frac{ff'}{r}\right]y^ty^t+\left[-\frac{f''}{2f}-\frac{f'}{rf}\right]y^r y^r+\left[\bar{R}ic-f-rf'\right]\bar{F}^2~
\end{eqnarray}
where $\bar{R}ic$ denotes the Ricci scalar of Finsler structure $\bar{F}$. We have used the property of homogenous function $H(\lambda y)=\lambda^n H(y)$, i.e. $y^\mu \frac{\partial H(y)}{\partial y^\mu}=nH(y)$, to derive the geodesic spray coefficients $G^\mu$ and Ricci scalar.
By equation (\ref{Ricci scalar1}), we obtain the solution of $Ric=0$. It is of the form
\begin{eqnarray}
\bar{R}ic&=&k,\\
f_S&=&k-2GM/r,
\end{eqnarray}
where $k=\pm1,0$. And the solution of constant Ricci scalar ($Ric=3b$) is given as
\begin{eqnarray}
\bar{R}ic&=&k,\\
\label{Schwar-dS}
f_{Sd}&=&k-2GM/r-br^2.
\end{eqnarray}

A Finsler spacetime with constant flag curvature $K$ must have constant Ricci scalar $(n-1)K$ ($n$ is the dimension of the spacetime). However, the reverse statement is not true. Now, we test whether the solution (\ref{Schwar-dS}) is corresponded to the Finsler spacetime with constant flag curvature or not. A Finsler spacetime with constant flag curvature $K$ is equivallent to its predecessor of flag curvature possess the following form \cite{Book by Bao}
\begin{equation}\label{const curv}
F^2R^\mu_\nu=K\left(F^2\delta^\mu_\nu-\frac{y^\mu}{2} \frac{\partial F^2}{\partial y^\nu}\right),
\end{equation}
where $F^2R^\mu_\nu$ is defined as
\begin{equation}\label{Ricci scalar2}
F^2R^\mu_{~\nu}=2\frac{\partial G^\mu}{\partial x^\nu}-y^\lambda\frac{\partial^2 G^\mu}{\partial x^\lambda\partial y^\nu}+2G^\lambda\frac{\partial^2 G^\mu}{\partial y^\lambda\partial y^\nu}-\frac{\partial G^\mu}{\partial y^\lambda}\frac{\partial G^\lambda}{\partial y^\nu}.
\end{equation}
Plugging the ansazt into the formula (\ref{Ricci scalar2}), after tedious calculation, we obtain
\begin{eqnarray}\label{Rtt}
F^2R^t_t&=&-\frac{f''}{2f}y^ry^r-\frac{rf'}{2}\bar{F}^2,\\
\label{Rtr}
F^2R^t_r&=&\frac{f''}{2f}y^ty^r,\\
\label{Rti}
F^2R^t_i&=&\frac{rf'}{4}y^t\frac{\partial\bar{F}^2}{\partial y^i},\\
\label{Rrr}
F^2R^r_r&=&\frac{ff''}{2}y^ty^t-\frac{rf'}{2}\bar{F}^2,\\
\label{Rri}
F^2R^r_i&=&\frac{rf'}{4}y^r\frac{\partial\bar{F}^2}{\partial y^i},\\
\label{Rij}
F^2R^i_j&=&\bar{F}^2\bar{R}^i_j+\left(\frac{ff'}{2r}y^ty^t-\frac{f'}{2rf}y^ry^r-f\bar{F}^2\right)\delta^i_j+\frac{f}{2}y^i\frac{\partial\bar{F}^2}{\partial y^j}.
\end{eqnarray}
Plugging the solution (\ref{Schwar-dS}) into the above formula (\ref{Rtt},\ref{Rtr},\ref{Rti},\ref{Rrr},\ref{Rri},\ref{Rij}), it is obvious from the formula (\ref{const curv}) that Finsler metric (\ref{Schwarz like}) with the solution (\ref{Schwar-dS}) does not have constant flag curvature if the parameters $M\neq0$ and $b\neq0$. However, it is Einstein metric since its Ricci scalar is constant. By making use of the property $R_{\mu\nu}=R_{\nu\mu}$ and noticing the $\bar{F}$ is a Finsler structure with constant Ricci scalar, we obtain that the special case of the solution (\ref{Schwar-dS}) with parameter $M=0$ and $b\neq0$ is corresponded to Finsler spacetime with constant flag curvature, i.e., it satisfies the following relation
\begin{equation}
F^2R^\mu_\nu=b\left(F^2\delta^\mu_\nu-\frac{y^\mu}{2} \frac{\partial F^2}{\partial y^\nu}\right).
\end{equation}
It means that the flag curvature is constant and equals to $b$.

Pfeifer et al.have studied the electromagnetic field in Finsler spacetime \cite{Pfeifer EM field}. However, no specific solution of electromagnetic field equation, such as static electric field, is discussed. Now, we study the solution of Finslerian gravity with electric charge. Analogy to the Riemannian Ressiner-Nordstrom metric, we suggest that the Finsler metric has the form of the ansatz (\ref{Schwarz like}) with $f_{RN}=k-\frac{2GM}{r}+\frac{4\pi_F GQ^2}{r^2}$. Then, by formula (\ref{Ricci scalar1}), its Ricci scalar is of the form
\begin{eqnarray}\label{Ricci scalar em}
F^2Ric=\frac{4\pi_F GQ^2}{r^4}\left(fy^ty^t-f^{-1}y^r y^r+r^2\bar{F}^2\right).
\end{eqnarray}
The gravitational field equation in the given Finsler spacetime (\ref{Schwarz like}) should be of the form \cite{Finsler BH}
\begin{equation}\label{field equation}
G^\mu_\nu=8\pi_F G T^\mu_\nu,
\end{equation}
where the modified Einstein tensor $G_{\mu\nu}$ is defined as
\begin{equation}\label{Einstein tensor}
G_{\mu\nu}\equiv Ric_{\mu\nu}-\frac{1}{2}g_{\mu\nu}S,
\end{equation}
and $4\pi_F$ denotes the volume of $\bar{F}$.
Here, the Ricci tensor we used is first introduced by Akbar-Zadeh\cite{Akbar}
\begin{equation}\label{Ricci tensor}
Ric_{\mu\nu}=\frac{\partial^2\left(\frac{1}{2}F^2 Ric\right)}{\partial y^\mu\partial y^\nu}.
\end{equation}
And the scalar curvature in Finsler geometry is given as $S=g^{\mu\nu}Ric_{\mu\nu}$. Then, by making use of the equations (\ref{Ricci scalar em},\ref{Einstein tensor},\ref{Ricci tensor}), we obtain the non vanishing components of Einstein tensor
\begin{equation}\label{Einstein tensor res}
Ric_{\mu\nu}=G_{\mu\nu}=\frac{4\pi_F GQ^2}{r^4}\rm diag\{f,-f^{-1},r^2\bar{g}_{ij}\}.
\end{equation}
One should notice that the scalar curvature $S$ vanishes. It means that the trace of energy momentum tensor vanishes. And such fact implies that the particle is massless in physics. Plugging the result of Einstein tensor (\ref{Einstein tensor res}) into the field equation (\ref{field equation}), we obtain the energy momentum tensor
\begin{equation}\label{electric EM tensor}
T_{\mu\nu}=\frac{Q^2}{2r^4}\rm diag\{f,-f^{-1},r^2\bar{g}_{ij}\}.
\end{equation}

In ref. \cite{Finsler BH}, we have proved that the covariant derivative of the modified Einstein tensor is conserved in the Finsler spacetime (\ref{Schwarz like}), i.e., $G^\mu_{\nu|\mu}=0$, where ``$|$" denotes the covariant derivative. Since the Finslerian Schwarzschild-de Sitter solution $f=f_{Sd}$ and the Finslerian Ressiner-Nordstrom solution $f=f_{RN}$ both possess the same form with Finslerian Schwarzschild solution (\ref{Schwarz like}). Thus, following the same process of ref. \cite{Finsler BH}, one can find that the modified Einstein tensor and the energy momentum tensor corresponded to the solution are also conserved. Energy momentum tensor of general electromagnetic field in Finsler spacetime should reduce to the energy momentum tensor given above (\ref{electric EM tensor}) if the Finsler spacetime reduce to Finslerian Ressiner-Nordstrom solution, i.e., the ansatz metric (\ref{Schwarz like}) with $f=f_{RN}$.

\subsection{Symmetry of Finslerian Ressiner-Nordstrom solution}
In Riemannian geometry, spaces with constant sectional curvature are equivalent. And all these spaces have $n(n+1)/2$ independent Killing vectors. However, The number of independent Killing vectors of an $n$ dimensional non-Riemannian Finsler spacetime should be no more than $\frac{n(n-1)}{2}+1,~ n\neq4 ~\&~ n\geq3$ \cite{Wang}. And a four dimensional non-Riemannian Finsler spacetime has no more than 8 independent Killing vectors \cite{Bogoslovsky}. In general, Finslerian extension of a given Riemannian spacetime has lower symmetry than the Riemannian spacetime. For example, the Finslerian Schwarzschild spacetime $f=f_S$ only has two independent Killing vectors\cite{Finsler BH}. The Kerr spacetime has two independent Killing vectors. Thus, Finslerian extension of Kerr spacetime should only has one or zero independent Killing vector. Furthermore, if we require the Finslerian extensional spacetime is static, then, arbitrary static Finsler spacetime should be Finslerian extension of Kerr spacetime. This fact implies that it is hard to find Finslerian extension of Kerr spacetime. In Ref.\cite{Finsler BH}, we have shown that the Finsler spacetime with specific form (\ref{Schwarz like}) could form a horizon at $f=0$. Therefore, at present, one can find from Finslerian Ressiner-Nordstrom solution that Finslerian black hole has two properties, namely, mass and electric charge.

The Killing equation $K_V(F)$ in Finsler spacetime is of the form\cite{Finsler PF}
\begin{equation}
\label{killing F}
K_V(F)\equiv \tilde{V}^\mu\frac{\partial F}{\partial x^\mu}+y^\nu\frac{\partial \tilde{V}^\mu}{\partial x^\nu}\frac{\partial F}{\partial y^\mu}=0.
\end{equation}
Plugging the formula (\ref{Finsler metric}) into the Killing equation (\ref{killing F}), we obtain that
\begin{equation}\label{killing g}
V^\mu\frac{\partial g_{\alpha\beta}}{\partial x^\mu}+g_{\alpha\lambda}\frac{\partial V^\lambda}{\partial x^\beta}+g_{\lambda\beta}\frac{\partial V^\lambda}{\partial x^\alpha}+y^\nu\frac{\partial V^\mu}{\partial x^\nu}\frac{\partial g_{\alpha\beta}}{\partial y^\mu}=0.
\end{equation}
The left side of Killing equation (\ref{killing g}) is just Lie derivative of Finsler metric $g_{\alpha\beta}$ \cite{Yano book}. From the equation (\ref{killing g}), we investigate the symmetry of the Finsler spacetime (\ref{Schwarz like}) with constant flag curvature, i.e., $f=f_d=k-br^2$. It should be noticed that the Killing equation (\ref{killing g}) differs from Riemannian Killing equation in $y^\nu\frac{\partial V^\mu}{\partial x^\nu}\frac{\partial g_{\alpha\beta}}{\partial y^\mu}$. And only component $g_{ij}$ of the Finsler metric of the Finsler spacetime (\ref{Schwarz like}) has $y$ dependence. Also, Finsler spacetime (\ref{Schwarz like}) with $f=f_d$ reduces to Riemannian spacetime with constant sectional curvature if $\bar{F}$ reduces to Riemannian surface with constant sectional curvature. Therefore, there are three independent Killing vectors which have index $t$ and $r$ only, and these Killing vectors only depend on coordinate $t$ and $r$. And we have shown in Ref.\cite{Finsler BH} that a ``Finslerian sphere" with constant positive flag curvature admits one Killing vector $V^\varphi=C^\varphi$, where $C^\varphi$ is a constant. It is obvious that $V^\varphi=C^\varphi$ is Killing vector of the Finsler spacetime (\ref{Schwarz like}) with constant flag curvature. Finally, we conclude that the Finsler spacetime
\begin{equation}\label{Finsler de Sitter}
F^2_d=-(1-br^2)y^ty^t+(1-br^2)^{-1}y^ry^r+r^2F^2_{FS}
\end{equation}
admits four independent Killing vectors. The specific form of ``Finslerian sphere" will be given in next section. We have shown in Ref.\cite{Finsler PF} that a four dimensional projectively flat Randers spacetime with constant flag curvature admits six independent Killing vectors. Thus, Finsler spacetime (\ref{Finsler de Sitter}) and projectively flat Randers spacetime with constant flag curvature are not equivalent, namely, after a coordinate transformation, one can change as other. This fact is quite different with the one in Riemannian geometry. Since each Riemannian spacetime with constant sectional curvature are equivalent.

\section{Finslerian Laplace operator on ``Finslerian sphere"}
Spherical harmonics is the eigenfunction of Laplace operator of Riemannian 2-sphere. Since our Finslerian Ressiner-Nordstrom solution in Finsler spacetime differs from Ressiner-Nordstrom metric only in $\bar{F}$. It is worth to investigate the Laplace operator of Finslerian surface with positive constant flag curvature. In the discussion of above section, we have shown that Finsler spaces with constant flag curvature may not equivalent to each other. Thus, we adopt a specific form of Finsler surface with positive constant flag curvature to study its Laplace operator, i.e., a two dimensional Randers-Finsler space with constant positive flag curvature \cite{Bao}
\begin{equation}\label{Finsler sphere}
F_{\rm FS}=\frac{\sqrt{(1-\epsilon^2\sin^2\theta)y^\theta y^\theta+\sin^2\theta y^\varphi y^\varphi}}{1-\epsilon^2\sin^2\theta}-\frac{\epsilon\sin^2\theta y^\varphi}{1-\epsilon^2\sin^2\theta},
\end{equation}
where $0\leq\epsilon<1$. We call it ``Finslerian sphere".

In Riemannian geometry, the Laplace operator can be defined in several different ways \cite{Gallot} and these Laplace operators are equivalent to each other. However, Finslerian extension of these definitions will lead to different Laplace operators. Various Finslerian Laplace operators are respectively defined by Bao and Lackey\cite{Bao and Lackey}, Shen\cite{Shen's book}, Barthelme \cite{Barthelme}. In this paper, we will adopt the Finslerian Laplace operator defined by Barthelme. The Finslerian Laplace operator for ``Finslerian sphere" (\ref{Finsler sphere}) is of the form \cite{Barthelme}
\begin{eqnarray}
\Delta_{FS}&=&\frac{2(1-\epsilon^2\sin^2\theta)^{3/2}}{\sin^2\theta(1+\sqrt{1-\epsilon^2\sin^2\theta})}\frac{\partial^2}{\partial\varphi^2}+\frac{2(1-\epsilon^2\sin^2\theta)}{1+\sqrt{1-\epsilon^2\sin^2\theta}}\frac{\partial^2}{\partial\theta^2}\nonumber\\
\label{Finsler spherical harmonic}
&&+\frac{2\cos\theta(\epsilon^2\sin^2\theta+\sqrt{1-\epsilon^2\sin^2\theta})}{\sin\theta(1+\sqrt{1-\epsilon^2\sin^2\theta})}\frac{\partial}{\partial\theta}.
\end{eqnarray}
While $\epsilon=0$, the Finslerian Laplace operator (\ref{Finsler spherical harmonic}) reduces to Riemanian Laplace operator for 2-sphere. And its eigenfunction is just spherical harmonics $Y_{lm}(\theta,\varphi)$ which satisfies
\begin{equation}\label{spherical harmonics}
\Delta_{FS}|_{\epsilon=0}Y_l^m=-l(l+1)Y_l^m.
\end{equation}
By using the equation (\ref{spherical harmonics}), one can find that $Y_0^0,Y_1^0,Y_1^{\pm1}$ are eigenfunctions of the Finslerian Laplace operator $\Delta_{FS}$, and they satisfy the following equations
\begin{eqnarray}
\Delta_{FS}Y_0^0=0,~\Delta_{FS}Y_1^0=-2Y_1^0,~\Delta_{FS}Y_1^{\pm1}=(-2+2\epsilon^2)Y_1^{\pm1}~.
\end{eqnarray}
It is obvious that the eigenvalues corresponded to $Y_0^0,Y_1^0$ remain the same with Riemannian Laplace operator and the eigenvalues corresponded to $Y_1^{\pm1}$ are different. This fact just reflect the symmetry of ``Finslerian sphere". As for the symmetry of ``Finslerian sphere", We have showed that $z$-axis rotational symmetry, which corresponds to Killing vector $V^\varphi=C^\varphi$, is preserved and other symmetry of Riemanian 2-sphere is broken \cite{Finsler BH}.

The symmetry of ``Finslerian sphere" may account for some specific physical phenomena, such as the power asymmetry of CMB \cite{Power asymmetry}. And such phenomena can be treated as a perturbation of standard physical theories. Thus, we expand Finslerian Laplace operator $\Delta_{FS}$ in powers of $\epsilon$. To first order in $\epsilon^2$, the Finslerian Laplace operator (\ref{Finsler spherical harmonic}) is given as
\begin{eqnarray}\label{Finsler spherical harmonic1}
\Delta_{FS}&=&\frac{4-5\epsilon^2\sin^2\theta}{4\sin^2\theta}\frac{\partial^2}{\partial\varphi^2}+\left(1-\frac{3}{4}\epsilon^2\sin^2\theta\right)\frac{\partial^2}{\partial\theta^2}
+\frac{\cos\theta}{\sin\theta}\left(1+\frac{3}{4}\epsilon^2\sin^2\theta\right)\frac{\partial}{\partial\theta}~.
\end{eqnarray}
By making use of the recurrence formula of spherical harmonics, to first order in $\epsilon^2$, we obtain the eigenfunction of the Finslerian Laplace operator (\ref{Finsler spherical harmonic1})
\begin{eqnarray}\label{eigenfunction}
\bar{Y}_l^m=Y_l^m+\epsilon^2(C_{l+2}^m Y_{l+2}^m+C_{l-2}^m Y_{l-2}^m),
\end{eqnarray}
where
\begin{eqnarray}
C_{l+2}^m&=&-\frac{3l(l-1)}{8(2l+3)^2}\sqrt{\frac{(l+m+1)(l-m+1)(l+m+2)(l-m+2)}{(2l+1)(2l+5)}}~,\\
C_{l-2}^m&=&\frac{3(l+1)(l+2)}{8(2l-1)^2}\sqrt{\frac{(l+m)(l-m)(l+m-1)(l-m-1)}{(2l+1)(2l-3)}}~.
\end{eqnarray}
And the corresponded eigenvalue of the Finslerian Laplace operator is given as
\begin{equation}\label{eigenvalue}
\lambda=-l(l+1)+\epsilon^2\left(\frac{3(l-1)l(l+1)(l+2)}{2(2l-1)(2l+3)}+\frac{m^2(14l^3+21l^2+19l+6)}{2(2l+1)(2l-1)(2l+3)}\right)~.
\end{equation}

The spatial geometry of the Finslerian Schwarzschild metric is
\begin{equation}\label{3d metric}
F^2_{3d}=y^ry^r+r^2\bar{F}^2_{FS},
\end{equation}
while $M=0$. Since $\bar{F}_{FS}$ does not depend on $r$, thus, following the definition of Finslerian Laplace operator \cite{Barthelme}, we find that Finslerian Laplace equation for the three dimensional space (\ref{3d metric}) is of the form
\begin{equation}\label{3d Laplace eq}
\Delta W=\frac{1}{r^2}\frac{\partial}{\partial r}\left(r^2\frac{\partial W}{\partial r}\right)+\Delta_{FS}W=0~.
\end{equation}
Solution of the Finslerian Laplace equation (\ref{3d Laplace eq}) is of the form
\begin{equation}\label{3d solution}
W=(Ar^{n_1}+Br^{n_2})\bar{Y}_l^m,
\end{equation}
where $A$ and $B$ are constant that depend on boundary condition of the Finslerian Laplace equation, and
\begin{eqnarray}
n_1&=&\frac{-1+\sqrt{1-4\lambda}}{2},\\
n_2&=&\frac{-1-\sqrt{1-4\lambda}}{2}.
\end{eqnarray}
Two facts for the solution (\ref{3d solution}) should be noticed. One is the index $n_1$ and $n_2$ depend not only on $l$ but also on $m$. Another fact is monopolar and dipolar term of multipole expansion of gravitational potential is unchanged in a Finsler spacetime with ``Finslerian sphere". It is consistent with our previous research given in ref.\cite{Finsler BH} where we have shown that the gravitational potential for the Schwarzschild-like spacetime is the same with Newtonian gravity in weak field approximation.
However, due to the Finslerian modification of the eigenvalue $\lambda$ (\ref{eigenvalue}), other multipole terms is changed.

\section{Discussions and conclusions}\label{sec:conclusions}
Beside the definition of Ricci tensor introduced by Akbar-Zadeh, Shen introduce another definition \cite{Shen Ric tensor}
\begin{equation}
Ric_{\mu\nu}=(R^{~\lambda}_{\mu~\lambda\nu}+R^{~\lambda}_{\nu~\lambda\mu})/2,
\end{equation}
where $R^{~\lambda}_{\mu~\lambda\nu}$ denotes the Riemann curvature tensor of the Berwald connection. The two definitions are equivalent if the Finsler spacetime has constant flag curvature. However, one can check that the two definitions are not equivalent for the ansazt metric with $f_S=k-\frac{2GM}{r}$, i.e., the Ricci flat case. Since the Ricci tensor introduced by Akbar-Zadeh corresponded to Finslerian gravitational field equation (\ref{field equation}) with exact solutions, such as $f=f_S$, $f=f_{Sd}$ and $f=f_{RN}$, therefore, the Ricci tensor introduced by Akbar-Zadeh is preferred in physics.

In Finsler geometry, there are two types of volume form, namely, the Busemann-Hausdorff volume form and Holmes-Thompson volume form \cite{Shen's book}. We have shown in Ref.\cite{Finsler BH} that the volume of ``Finslerian sphere" in terms of the Busemann-Hausdorff volume form is $4\pi$. The Finslerian Laplace operator for ``Finslerian sphere" (\ref{Finsler spherical harmonic}) is defined on fiber bundle and its definition is related to Holmes-Thompson volume form \cite{Barthelme}. The volume of ``Finslerian sphere" in terms of Holmes-Thompson volume form is given as
\begin{equation}
{\rm Vol}_{FS}=\int\frac{\sin\theta}{(1-\epsilon^2\sin^2\theta)^{3/2}}d\theta\wedge d\varphi=\frac{4\pi}{1-\epsilon^2}.
\end{equation}
The two definitions of volume form will slightly alter the Finslerian gravitation field equation (\ref{field equation}), for the term $8\pi_F$ is double the surface volume of Finsler space $\bar{F}$. In general relativity, black hole entropy depends on the surface volume of the black hole. Thus, studying the black hole thermodynamics in Finsler spacetime will help us to find which volume form is preferred in physics. It will be discussed in our future work.

In this paper, we have obtained Finslerian Ressiner-Nordstrom solution where it is asymptotic to a Finsler spacetime with constant flag curvature while $r\rightarrow\infty$ (\ref{Finsler de Sitter}). The covariant derivative of modified Einstein tensor in Finslerian gravitational field equation for this solution is conserved. The symmetry of the special Finslerian Ressiner-Nordstrom spacetime (\ref{Finsler de Sitter}) has been investigated. It admits four independent Killing vectors. The Finslerian Ressiner-Nordstrom solution differs from Ressiner-Nordstrom metric only in two dimensional subspace $\bar{F}$. The Finslerian Ressiner-Nordstrom could form horizons at $f_{RN}=0$. At present, we can conclude that Finslerian black hole have at least two properties, namely, mass and electric charge. Our solutions show that two dimensional subspace $\bar{F}$ has constant flag curvature. Spherical harmonics is the eigenfunction of Laplace operator of Riemannian 2-sphere. We have obtained the eigenfunction of Finslerian Laplace operator introduced by Barthelme\cite{Barthelme} of ``Finslerian sphere" (\ref{Finsler sphere}). The eigenfunction (\ref{eigenfunction}) is just a combination of spherical harmonics in powers of Finslerian parameter $\epsilon$. However, the eigenvalue depends on both $l$ and $m$. The eigenvalues corresponded to $Y_1^0$ remain the same with Riemannian Laplace operator and the eigenvalues corresponded to $Y_1^{\pm1}$ are different. This fact just reflect the symmetry of ``Finslerian sphere", which admits a Killing vector $V^\varphi=C^\varphi$ and breaks other symmetry of Riemannian sphere. The eigenfunction (\ref{eigenfunction}) of Finslerian Laplace operator implies that monopolar and dipolar terms of multipole expansion of gravitational potential are unchanged and other multipole terms are changed.

\begin{acknowledgments}
This work has been supported by the National Natural Science Fund of China under grant Nos. 11775038.
\end{acknowledgments}

\end{document}